\magnification=1000
\hsize=11.7cm
\vsize=18.9cm
\lineskip2pt \lineskiplimit2pt
\nopagenumbers

\hoffset=-1truein
\voffset=-1truein

\advance\voffset by 4truecm
\advance\hoffset by 4.5truecm

\newif\ifenteteCH
\newif\ifenteteIN

\headline{\ifenteteCH\ifodd	\count0 
      \rlap{\headCH}\hfill\tenrm\llap{\the\count0}\relax
    \else
        \tenrm\rlap{\the\count0}\hfill\llap{\headIN} \relax
    \fi\else
\global\enteteCHtrue\fi}

\def\enteteCH#1{\enteteINfalse\gdef\headCH{#1}}
\def\enteteIN#1{\enteteCHfalse\gdef\headIN{#1}}
\enteteCH{}
\enteteIN{}

\input amssym.def
\input amssym.tex

\def\-{\hbox{-}}
\def\.{{\cdot}}

\def\F{{\cal F}}

\def\s{{\cal S}}

\def\Ab{\frak A\frak b}

\def\int{\frak i\frak n\frak t}

\def\qq{\quad{\rm and}\quad}

\def\too{\longrightarrow}

 3
 2
\font\large=cmr10  scaled \magstep 2
 2
\font\larti=cmti10  scaled \magstep 2
 2
\font\cds=cmr7

\centerline{\large Frobenius {\larti P}-categories via the Alperin condition}
\bigskip
\centerline{\bf Lluis Puig }

\smallskip
\centerline{\cds CNRS, Institut de Math\'ematiques de Jussieu}
\smallskip
\centerline{\cds 6 Av Bizet, 94340 Joinville-le-Pont, France}
\smallskip
\centerline{\cds puig@math.jussieu.fr}

\bigskip
\noindent
{\bf  £1. Introduction}

\bigskip
£1.1. Let $P$ be a finite $p\-$group and $\F$ a {\it divisible $P\-$category\/}. In [5,~Ch.~5] we showed that our approach 
 in [4, Appendix] to Alperin's Fusion Theorem~[1] for {\it local pointed groups\/}  can be translated to $\F$ and that, in this context, 
it still makes sense to define the {\it $\F\-$essential\/} subgroups of $P$ [5,~5.7]. Then, we rose the following question:  
{\it in what extend the  behaviour of the $\F\-$essential subgroups of $P$ characterizes the Frobenius 
$P\-$categories?\/} In this Note we give a more satisfactory answer to this question than what we obtained in~[5,~Theorem~5.22].

 \medskip
£1.2. Let us recall our setting. A {\it divisible $P\-$category\/} is a subcategory~$\F$ of
the category of finite groups containing the {\it Frobenius category~$\F_{\!P}$\/} of~$P$ [5,~1.8] --- where the  
objects are  all the subgroups of $P$  and where  all the homomorphisms are injective --- and fulfilling the following condition:
\smallskip
\noindent
£1.2.1.\quad {\it If $Q\,,$ $R$ and $T$ are subgroups of $P\,,$ for any
$\varphi\in \F (Q,R)$ and any group homomorphism $\psi\,\colon T\to R$
the composition $\varphi\circ\psi$ belongs to $\F (Q,T)$ (if~and) only if
$\psi\in \F (R,T)\,.$\/}
\smallskip
\noindent
Here, $\F (Q,R)$ denotes the set of $\F\-$morphisms from $R$ to $Q\,.$ Moreover,  we  consider the category
${\Bbb Z}\F$ still defined over the set of all the subgroups of~$P$ where, for any pair of
subgroups  $Q$ and~$R$ of~$P\,,$ the set of morphisms from~$R$ to
$Q$ is the {\it free ${\Bbb Z}\-$module\/} ${\Bbb Z}\F (Q, R)$ over~$\F (Q,R)\,,$  with the
{\it distributive\/} composition extending the composition in~$\F\,.$

\medskip
£1.3. For any two different elements $\varphi\,,\varphi' \in \F (Q,R)\,,$
we call $\F\-${\it dimor-phism\/} from $R$ to $Q$ the difference $\varphi' - \varphi$ in
the $\Bbb Z\-$module $\Bbb Z\F (Q,R)\,;$ it is clear that  the set of 
$\F\-$dimorphisms is stable by left-hand and right-hand composition with $\F\-$morphisms;
note that, for any $\varphi\in \F (Q,R)\,,$ the family 
$\{\varphi'-\varphi\}_{\varphi'} \,,$ where $\varphi'$ runs over $\F (Q, R)-\{\varphi\}\,,$
is a ${\Bbb Z}\-$basis of the kernel of the evident {\it augmentation\/}
${\Bbb Z}\-$linear map
$$\varepsilon_{Q,R} : {\Bbb Z} \F (Q,R)\too {\Bbb Z}
\eqno £1.3.1\phantom{.}$$ 
sending any $\varphi\in \F (Q,R)$ to $1\,.$

\medskip
£1.4. The next elementary lemma  [5, Lemma~5.4] relates any ``linear'' decomposition of an 
$\F\-$dimorphism in terms of $\F\-$dimorphisms with the old {\it partially defined  linear combinations\/} introduced in~[3,~Ch.~III]. 
Note that, in the case where $Q = P\,,$ $\varphi$ is the inclusion map $\iota_R^P\,\colon R\to P\,,$ 
and for any $i\in I\,,$ we have $\mu_i = \iota_{Q_i}^P\,,$  $Q_i = R_i$ and $\varphi_i = {\rm id}_{R_i}\,,$ 
equalities~£1.5.2 below coincide  with the decomposition pattern in the original formulation of Alperin's Fusion Theorem~[1].

\bigskip
\noindent
{\bf Lemma £1.5.} {\it With the notation above, let $\{Q_i\} _{i\in I}$
and $\{R_i\}_{i\in I}$ be finite fami-lies of subgroups of $P$ and, for
any $i\in I\,,$ let $\varphi'_i -\varphi_i$ be an $\F\-$dimorphism from $R_i$ to $Q_i$ and 
$\mu_i\,\colon Q_i\to Q$ and  $\nu_i\,\colon R\to R_i$ be two $\F\-$morphisms. Then, we have
$$\varphi' -\varphi = \sum_{i\in I} \mu_i\circ
(\varphi'_i -\varphi_i)\circ \nu_i
\eqno £1.5.1\phantom{.}$$
if and only if there are $n\in {\Bbb N}$ and an injective map $\sigma\,\colon
\Delta_n\to I$ fulfilling 
$$\eqalign{\varphi &= \mu_{\sigma (0)}\circ\varphi_ {\sigma (0)}\circ \nu_{\sigma (0)}\cr
\mu_{\sigma(\ell -1)}\circ \varphi'_{\sigma(\ell -1)}\circ\nu_{\sigma(\ell -1)}  
&= \mu_{\sigma (\ell)}\circ\varphi_ {\sigma (\ell)} \circ
\nu_{\sigma (\ell)}\;\; \hbox{for any $1\le \ell\le n$}\cr 
\mu_{\sigma (n)}\circ \varphi'_{\sigma (n)}\circ
\nu_{\sigma (n)}   &= \varphi'\,.\cr}
\eqno £1.5.2\phantom{.}$$ \/}

\medskip
£1.6. According to {\it Yoneda's Lemma\/} [2,~\S1],  the {\it contravariant\/} functor $\frak h_\F\,\colon \F \to \Ab$ mapping any
subgroup $Q$ of $P$ on $\Bbb Z\F (P,Q)$ and any $\F\-$morphism $\varphi\,\colon R\to Q$ on the
group homomorphism $\frak h_\F (Q)\to \frak h_\F (R)$ defined by the composition with
$\varphi$ is a {\it projective object\/} in the category of {\it contravariant\/} functors from $\F$ to $\Ab\,.$
Then, denoting by $\Bbb Z\,\colon \F\to \Ab$ the trivial {\it contravariant\/}
functor mapping any $\F\-$object on $\Bbb Z\,,$ the ring of integers, and any $\F\-$morphism on ${\rm
id}_\Bbb Z\,,$ the family of {\it augmentation maps\/} $\varepsilon_{P,Q}$ when $Q$ runs over
the set of subgroups of $P$ defines a surjective {\it natural map\/}
$$\varepsilon_\F : \frak h_\F\too \Bbb Z
\eqno £1.6.1\phantom{.}$$
and we set $\frak w_\F = {\rm Ker}(\varepsilon_\F)\,,$ which is nothing but the {\it Heller
translated\/} of the trivial functor $\Bbb Z\,.$

\medskip
£1.7. On the other hand, if $\frak a\,\colon \F\to \Ab$ is a {\it contravariant\/} functor,
 let us say that a family $\s = \{S_Q\}_Q$
of subsets $S_Q\i \frak a (Q)\,,$ where $Q$ runs over the set of proper subgroups of $P\,,$
is a {\it generator family of  $\frak a$\/} whenever, for any proper subgroup~$Q$
of $P\,,$ we have
$$\frak a (Q) = \sum_R\,\sum_{\varphi\in \F (R,Q)}\,\sum_{a\in S_R} \Bbb Z\.\big(\frak
a(\varphi)\big)(a)
\eqno £1.7.1,$$
 where $R$ runs over the set of subgroups of $P$ (such that $\vert R\vert\ge 
\vert Q\vert\,$). Now, it is quite clear from Lemma~£1.5 above that Alperin's Fusion Theorem~[1] provides
a particular {\it generator family\/} of the {\it Heller translated\/} $\frak w_\F\,.$

\medskip
£1.8.  In order to find minimal {\it generator families\/} of $\frak w_\F\,,$ let us define 
 a subfunctor $\frak r_\F$ of $\frak w_\F$ mapping any subgroup $Q$ of~$P$~on 
$$\frak r_\F (Q) = w_\F (P)\circ \Bbb Z\F(P,Q) +\sum_R \frak w_\F (R)\circ {\Bbb Z} \F (R,Q) 
\eqno £1.8.1\phantom{.}$$
\eject
\noindent
where $R$ runs over the set of subgroups of $P$ such that $\vert R\vert > \vert Q\vert\,.$ Then, 
we say that $Q$ is  $\F\-${\it essential\/} whenever~$\frak r_\F (Q)\not= \frak w_\F (Q)$ and
call $\F\-${\it irreducible\/} the elements of  $\frak w_\F (Q) - \frak r_\F (Q)\,.$ Coherently, the elements 
of $\frak r_\F (Q)$ are called  $\F\-${\it reducible\/}; actually, any element of $\frak r_\F (Q)$ is a sum of a family of
$\F\-$reducible  {\it $\F\-$dimorphisms\/} from $Q$ to $P\,.$ The following result [5,~Proposition~5.9]
justifies all these definitions.

\bigskip
\noindent
{\bf Proposition £1.9.}  {\it Let $\s =\{S_Q\}_Q$ be a generator family of\/
$\frak w_\F\,,$ where $Q$ runs over the set of proper subgroups of $P\,.$
The family formed by the $\F\-$irreducible elements of $S_Q\,,$ where $Q$ runs over the
set of $\F\-$essential subgroups of $P\,,$ is also a generator family of  
$\,\frak w_\F\,.$ Moreover, for any $\F\-$essential subgroup~$Q$ of~$P\,,$ there is
$\,\varphi\in\F (P,Q)$ such that\/ $S_{\varphi (Q)}$ contains an\/ $\F\-$irre-ducible element.\/}

\medskip
£1.10. Before going further, recall that $\F$ is a {\it Frobenius $P\-$category\/} whenever it fulfills 
the following two conditions [5,~2.8]
\smallskip
\noindent
{\bf Sylow condition.} {\it The group $\F_{\! P}(P)$ of inner automorphisms of
$P$ is a Sylow $p\-$subgroup of $\F (P)\,.$\/}
\smallskip
\noindent
{\bf Extension condition.} {\it For any subgroup $Q$ of $P\,,$ any subgroup $K$ of
${\rm Aut} (Q)$ and any $\F\-$morphism $\varphi\,\colon Q\to P$ such that $\varphi
(Q)$ is fully ${{}^{\varphi}\!}K\-$normalized in~$\F\,,$~there are an
$\F\-$morphism  $\psi\,\colon Q\.N_P^K (Q)\to P$ and $\chi\in K$ such that
$\psi (u) = \varphi \big(\chi (u)\big)$ for any $u\in Q\,.$\/}
\smallskip
\noindent
Here we set $\F(Q) = \F(Q,Q)\,,$ ${{}^{\varphi}\!}K$ denotes the image of $K$ in ${\rm Aut}\big(\varphi(Q)\big)$ throughout
the isomorphism $Q\cong \varphi(Q)$ induced by $\varphi$ and we say that {\it $Q$ is fully $K\-$normalized
in $\F$\/} whenever it fulfills [5,~2.6]
\smallskip
\noindent
£1.10.1\quad {\it For any $\psi\in \F \big(P,Q\.N_P^K (Q)\big)$ , we
have $\psi\big(N_P^K (Q)\big) = N_P^{{{}^\psi\!}K} \big(\psi (Q)\big)\,.$\/}
\smallskip
\noindent
Recall that we  say {\it fully centralized\/} or {\it fully normalized\/} whenever $K =\{1\}$ or~$K = {\rm Aut}(Q)\,,$
replacing $N^K_\F$ by $C_\F$ or $N_\F\,.$

\medskip
£1.11. According to Proposition~£1.9, when considering the {\it generator fami-lies\/} of $\frak w_\F\,,$
it suffices to consider the $\F\-$essential subgroups of $P\,.$ Now, if $Q$ is an $\F\-$essential subgroup of $P\,,$ 
we have $\frak h_\F(Q)/\frak r_\F(Q)\not\cong \Bbb Z$  and, denoting by~$\,\overline{\!\F (P,Q)\!}\,$ the image of $\F (P,Q)$ in the quotient $\frak h_\F(Q)/\frak r_\F(Q)\,,$ it is clear~that $\F (Q)$ acts on $\,\overline{\!\F(P,Q)\!}\,$ by composition on the left. At this point, it follows from [5,~Theorem~5.11] that:
\smallskip
\noindent
£1.11.1\quad {\it  If $\F$ is a Frobenius $P\-$category then $\F (Q)$ is transitive   on $\,\overline{\!\F (P,Q)\!}\,\,,$ $Q$ is {\it $\F\-$selfcentralizing\/}  and we have ${\bf O}_p\big(\F(Q)\big) = \F_Q(Q)\,.$\/}
\smallskip
\noindent
Recall that $Q$ is an {\it $\F\-$selfcentralizing\/} subgroup of~$P$ if $C_P\big(\varphi (Q)\big) 
= Z\big(\varphi (Q)\big)$ for any $\varphi \in \F (P,Q)$ [5,~4.8] and let us say that $Q$ is an {\it $\F\-$radical\/} if it is\break
\eject
\noindent 
$\F\-$selfcentralizing  and we have ${\bf O}_p\big(\F(Q)\big) = \F_Q(Q)\,.$  Moreover, recall that $Q$ is  an 
{\it $\F\-$intersected\/} subgroup of $P$  if it is selfcentralizing and fulfills [5,~4.11]
$$\F_Q (Q) = \bigcap_{\varphi \in \F (P,Q)} {}^{\varphi^*}\!\F_{\!P}\big(\varphi (Q)\big)
\eqno £1.11.2;$$
actually, an $\F\-$radical is an $\F\-$intersected subgroup of $P\,.$
 Note that statement~£1.11.1, Proposition~£1.9 and Lemma~£1.5 already prove the corresponding version in $\F$ of Alperin's 
 Fusion Theorem [5, Corollary~5.14]; thus, we consider 
the following condition on $\F\,:$
\smallskip
\noindent
{\bf Alperin condition.} {\it For any $\F\-$essential subgroup $Q$ of $P\,,$ $Q$ is an $\F\-$radical and $\F (Q)$ acts transitively    on $\,\overline{\!\F (P,Q)\!}\,\,.$\/}

\medskip
£1.12. On the other hand, for any subgroup $Q$ of $P$ and any subgroup $K$ of
${\rm Aut} (Q)$ such that $Q$ is fully $K\-$normalized in $\F\,,$ recall that the {\it divisible 
$N_P^K(Q)\-$subcategory $N_\F^K (Q)$ of $\F$\/} [5,~2.14] is the subcategory of $\F$ where, for any pair of subgroups~$R$
and~$T$ of~$N_P^K (Q)\,,$ the set of morphisms from $T$ to~$R$ is the set of elements
$\varphi\in\F (R,T)$ fulfilling the following condition:
\smallskip
\noindent
£1.12.1\quad {\it There are an $\F\-$morphism $\psi\, \,\colon Q\.T\to Q\.R$ and an element $\chi$ of $K$ 
such that $\chi (u) = \psi (u)$ for any $u\in Q$ and that $\psi (v) =\varphi (v)$ for any $v\in T\,.$\/}
\smallskip
\noindent
Note that, if $\F$ is a Frobenius $P\-$category then $N_\F^K (Q)$ is a Frobenius $N_\F^K (Q)\-$ category
too [5,~Proposition~2.16]. Our main purpose in this Note is to prove the following result. 

\bigskip
\noindent
{\bf Theorem~£1.13.} {\it A divisible $P\-$category is a Frobenius $P\-$category if and only if, for any
subgroup $Q$ of $P$ and any subgroup $K$ of ${\rm Aut} (Q)$ such that $Q$ is fully
$K\-$normalized in $\F\,,$ the $N_P^K(Q)\-$category $N_\F^K(Q)$ fulfills the Sylow and the Alperin conditions.\/}

\bigskip
\noindent
{\bf £2. Auxiliary results}

\bigskip
£2.1. In order to prove  Theorem~£1.13, it is handy to consider {\it partial Frobenius $P\-$categories\/} in the following sense. First of all, for short we  say that a triple $(Q,K,\varphi)$ formed by a subgroup~$Q$ of $P\,,$ a subgroup $K$ of ${\rm Aut} (Q)$ and an  $\F\-$morphism  $\varphi\,\colon Q\to P$ is {\it extensile\/} whenever there are an $\F\-$morphism  
$\psi\,\colon Q\.N_P^K (Q)\to P$ and an element $\chi$ of $K\cap \F (Q)$ such that $\psi (u) = \varphi \big(\chi (u)\big)$ 
for any $u\in Q\,;$ thus, the {\it extension condition\/} above states that  such a triple $(Q,K,\varphi)$ which fulfills that  $\varphi (Q)$ is fully ${{}^{\varphi}\!}K\-$normalized in~$\F$ is {\it extensile\/}.

\medskip
£2.2. Let $\frak X$ be a nonempty set of subgroups of $P$ containing any subgroup $Q$ of $P$ such that 
$\F(Q,R)\not= \emptyset$ for some $R\in \frak X\,,$ and denote by~$\F^{^\frak X}$ the {\it full\/} subcategory of 
$\F$ over $\frak X\,;$ we say that $\F^{^\frak X}$is a {\it partial Frobenius $P\-$category\/} if {\it $\F$ fulfills the Sylow condition and any triple  $(Q,K,\varphi)$ formed by an element $Q$ of $\frak X\,,$ a subgroup $K$ of ${\rm Aut} (Q)$ and an 
$\F\-$morphism  $\varphi\,\colon Q\to P$ such that $\varphi (Q)$ is fully ${{}^{\varphi}\!}K\-$normalized in~$\F$ is extensile\/}. From the proof of [5,~Corollary~2.13] it is straightforward to prove the following criterion that we need here.

\bigskip
\noindent
{\bf Proposition~£2.3.} {\it  With the notation above, assume that $\F$ fulfills the Sylow condition. Then, 
$\F^{^\frak X}$ is a partial Frobenius $P\-$category if and only if it fulfills the following two conditions
\smallskip
\noindent
{\rm £2.3.1}\quad For any pair of $\F\-$isomorphic elements $Q$ and $Q'$ of $\frak X\,,$ which are both
fully normalized and fully centralized in $\F\,,$ there is an $\F\-$isomorphism
$N_P (Q)\cong N_P (Q')$ mapping $Q$ onto~$Q'\,.$
\smallskip
\noindent
{\rm £2.3.2}\quad For any element $Q$ of $\frak X$ fully normalized and fully
centralized in $\F$ and any subgroup~$R$ of $N_P (Q)$ containing  $Q\.C_P (Q)\,,$
denoting by $\F (R)_Q$ the stabilizer of~$Q$ in~Ê$\F (R)\,,$ the  group
homomorphism $\F (R)_Q\to N_{\F (Q)}\big(\F_R (Q)\big)$
induced by the restriction is surjective.\/}

\medskip
£2.4. Similarly, note that all the definitions in~£1.6, £1.7 and~£1.8 above can be done in~$\F^{^\frak X}$ and then
an element $Q$ of $\frak X$ is $\F^{^\frak X}\-$essential if and only if it is $\F\-$essential; moreover, if 
$\F^{^\frak X}$ is a partial Frobenius $P\-$category, the characterization of the $\F\-$essential subgroups $Q$ in 
[5,~Theorem~5.11] still holds in~$\F^{^\frak X}\,.$  Here, we also need the lemma [5, Lemma~4.13] which can be restated as follows.

\bigskip
\noindent
{\bf Lemma £2.5.} {\it With the notation above, assume that $\F^{^\frak X}$ is a partial Frobenius $P\-$category.
Then, a triple $(R,J,\psi)$ formed by a subgroup $R$ of $P\,,$ a subgroup $J$ of ${\rm Aut} (R)$ and an $\F\-$morphism
$\psi\,\colon R\to P$ such that $\psi (R)$ is fully ${{}^{\psi}\!}J\-$normalized in $\F$
is extensile provided there are $Q\in \frak X$ having $R$ as a normal subgroup and stabilizing $J\,,$ and an 
$\F\-$morphism $\eta\,\colon Q\to P$ extending~$\psi\,.$\/}

\medskip
£2.6. Finally, we need the following characterization of the Frobenius $P\-$categories [5,~Theorem~4.12].

\bigskip
\noindent
{\bf Theorem £2.7}. {\it A divisible $P\-$category $\F$ fulfilling the Sylow condition is a Frobenius
$P\-$category  if and only if the following two conditions hold
\smallskip
\noindent
{\rm £2.7.1}\quad If $Q$ is an $\F\-$intersected subgroup of $P\,,$
$R$ is a subgroup of $N_P (Q)$ containing~$Q$ and $\varphi\,\colon Q\to P$ is an
$\F\-$morphism fulfilling  ${{}^{\varphi}}\F_{\! R} (Q)\i \F_{\! P}\big(\varphi
(Q)\big)$ then there is an $\F\-$morphism $\psi\,\colon R\to P$
extending~$\varphi\,.$
\smallskip
\noindent
{\rm £2.7.2}\quad Any divisible $P\-$category $\F'$ fulfilling $\F' (P,Q)\j\F (P,Q)$
for every $\F\-$intersected subgroup $Q$ of $P$ contains $\F\,.$\/}

\bigskip
\noindent
{\bf £3. Proof of Theorem~£1.13.}

\bigskip
 £3.1.  If $\F$ is a Frobenius $P\-$category then it follows from [5,~Proposition~2.16] that the $N_P^K(Q)\-$category $N_\F^K(Q)$ above is a Frobenius  $N_P^K(Q)\-$ca-tegory 
 and therefore it fulfills the Sylow and the Alperin conditions (cf.~£1.10 and~£1.11).

 \medskip 
£3.2.  Conversely, assume that for any subgroup $Q$ of $P$ and any subgroup~$K$ of ${\rm Aut} (Q)$
such that $Q$ is fully $K\-$normalized in~$\F\,,$ the $N_P^K(Q)\-$cate-gory $N_\F^K(Q)$ fulfills the Sylow and the Alperin conditions; we argue by induction on $\vert P\vert\,,$ $\prod_Q \vert \F (P,Q)\vert$ where 
$Q$ runs over the set of subgroups of $P\,,$ and~$\vert\frak X\vert$ successively;
since $\F$ fulfills the Sylow condition, we may assume that $\F^{^\frak X}$ is a partial Frobenius $P\-$category 
but $\frak X$ does not coincide with the set of all the subgroups of $P\,.$

\medskip
£3.3. Let $Q$ be a maximal subgroup of $P$ which does not belong to $\frak X\,;$  setting 
$$\frak Y = \frak X \cup  \{\varphi(Q)\}_{\varphi\in \F (P,Q)}
\eqno £3.3.1,$$
 it suffices to prove that $\F^{^\frak Y}$ fulfills both conditions in Proposition~£2.3 above;
actually, we may assume that $Q\not= \{1\}\,.$
 Let $\varphi\,\colon Q\to P$ be an $\F\-$morphism, set $Q' = \varphi (Q)$ and assume that
  $Q$ and $Q'$ are different and both fully normalized and fully centralized in $\F\,;$  then, according either to the very definition of $\F\-$essential subgroup or to the Alperin condition, in both cases the image $\bar\varphi$ of~$\varphi$ in $\,\overline{\!\F(P,Q)\!}\,$ coincides with 
$\,\overline{\!\iota_Q^P\circ\sigma}$ for some $\sigma\in \F (Q)\,,$ where $\iota_Q^P$ denotes the
inclusion map $Q\to P\,;$ that is to say, the difference $\varphi - \iota_Q^P\circ\sigma$ is
$\F\-$redu-cible and therefore we have (cf.~£1.8.1)
$$\varphi - \iota_Q^P\circ\sigma =  \sum_R\sum_{\theta\in  \frak w_\F
(R)}\,\theta\circ\alpha_{_{R,\theta}}
\eqno £3.3.2\phantom{.}$$
 for suitable $\alpha_{_{R,\theta}}\in \Bbb Z\F(R,Q)\,,$ where $R$ runs over the set of
subgroups of $P$ such that $\vert R\vert >\vert Q\vert\,.$

\smallskip
£3.4. Consequently, it follows from~£1.3 that we still have
$$\varphi - \iota_Q^P\circ\sigma =  \sum_{j\in J} (\psi'_j - \psi_j) \circ \mu_j 
\eqno £3.4.1\phantom{.}$$
where $J$ is a nonempty finite set and where, for any $j\in J\,,$ $\psi_j$ and $\psi'_j$ are elements of $\F (P,R_j)$
and $\mu_j\in \F (R_j,Q)$ for a suitable subgroup~$R_j$ of $P$ such that $\vert
R_j\vert >\vert Q\vert\,;$  more precisely, applying again the Alperin  condition and arguing
by induction on~$\vert P\,\colon Q\vert\,,$ we actually get
$$\varphi - \iota_Q^P\circ\sigma = \sum_{i\in I}\iota_{U_i}^P\circ  (\tau_i - 
{\rm id}_{U_i}) \circ \nu_i 
\eqno £3.4.2\phantom{.}$$
where $I$ is a nonempty finite set and, for any $i\in I\,,$ $\tau_i$ is an element of $\F (U_i)$ and
$\nu_i\in \F (U_i,Q)$ for a suitable subgroup~$U_i$ of~$P$ such that $\vert
U_i\vert >\vert Q\vert\,.$

\medskip
£3.5.  Then, it follows from Lemma~£1.5 that, for a suitable
$\,\ell\,,$ we can identify $\Delta_\ell$ with a subset of $I$ in such a way that $Q_0 = Q\,,$
$Q_{i +1}= \tau_i (Q_{i})$ for any $i\in \Delta_\ell\,,$ $Q_{\ell +1} = Q'$ and, denoting
by $\varphi_i\,\colon Q_{i}\cong Q_{i+1}$ the $\F\-$isomorphism induced by~$\tau_i\,,$ the
composition of all these isomorphisms coincides with the isomorphism $Q\cong Q'$ induced by
$\varphi\circ\sigma^{-1}\,.$ Moreover, note that $U_i$ contains $Q_i$ and $Q_{i+1}$ for
any $i\in \Delta_\ell$ and, in particular, it  belongs to~$\frak X\,.$
\eject

\medskip
£3.6. For any $i\in \Delta_{\ell +1}\,,$ let us  choose $\eta_i\in \F \big(P,N_P(Q_i)\big)$ such that $R_i =
\eta_i (Q_i)$ is fully normalized in $\F$ [5,~Proposition~2.7] and we may assume that~$R_0 = Q\,,$ that
$R_{\ell+1} = Q'$ and that $\eta_0$ and $\eta_{\ell+1}$ are the corresponding inclusion maps; moreover,
for any $i\in \Delta_\ell\,,$ denote by $\psi_i\,\colon R_{i}
\cong R_{i+1}$ the $\F\-$morphism mapping $\eta_{i}(u)$ on~$\eta_{i+1}\big(\varphi_{i}
(u)\big)$ for any $u\in Q_{i}\,.$ Then, for
any $i\in \Delta_\ell$  we claim that we can apply Lemma~£2.5 above to the triple $(R_{i},{\rm Aut} (R_{i}),\psi_i)\,;$ indeed, we are assuming that $\F^{^\frak X}$ is a partial Frobenius
$P\-$category; moreover, it is clear that  $R_{i}$ is a proper normal subgroup of
$\eta_{i}\big(N_{U_{i}}(Q_{i})\big)$ which clearly stabilizes ${\rm Aut}(R_{i})\,;$ 
finally, the $\F\-$morphism
$$\eta_{i}\big(N_{U_i}(Q_{i})\big)\too \eta_{i +1}\big(N_{U_i}(Q_{i+1})\big)
\eqno £3.6.1\phantom{.}$$
mapping $\eta_{i}(v)$ on $\eta_{i+1}\big(\tau_i (v)\big)$ for any $v\in N_{U_i}(Q_{i})$
clearly extends $\psi_i\,.$

\medskip
£3.7. Hence, since $R_{i+1}$ is fully normalized in $\F\,,$ it follows from this lemma that there
is an $\F\-$morphism $\zeta_{i}\,\colon N_P (R_{i})\to P$ extending
$\chi_i\circ\psi_i$ for some $\chi_i\in \F (R_i)\,;$ moreover, since $R_{i}$ is
fully normalized in $\F\,,$ we actually get
$$\zeta_{i}\big(N_P (R_{i})\big) = N_P (R_{i+1})
\eqno £3.7.1,$$
so that $\zeta_{i}$ induces an $\F\-$isomorphism $\xi_i\,\colon N_P (R_{i})\cong
N_P (R_{i+1})$ mapping $R_i$ onto $R_{i+1}\,.$ Finally, the composition of all these $\F\-$isomorphisms
when $i$ runs over  $\Delta_\ell$ yields an $\F\-$isomorphism $N_P (Q)\cong N_P (Q')$ which maps
$Q$ onto $Q'\,,$ proving condition~£2.3.1.

\medskip
£3.8. In order to prove condition~£2.3.2, we set $P' = N_P (Q)$ and we claim that the $P'\-$category $ \F' =N_\F (Q)$ still fulfills our hypothesis in~£3.2 above; more explicitly, if $R$ is a subgroup of $P'$ and~$J$ a subgroup
 of ${\rm Aut} (R)$  such that $R$ is fully $J\-$normalized in~$\F'\,,$ we claim 
 that the $N_{P'}^J (R)\-$category $N_{\F'}^J (R)$ fulfills the Sylow and the Alperin conditions. Set $T = Q\.R$
and denote by $I$  the subgroup of automorphisms of~$T$ which stabilize $Q$ and $R\,,$ and act on $R$ {\it via\/}
elements of $J\,;$ then, from its very definition (cf.~£1.12), it is easily checked that
$$N_P^I (T) = N_{P'}^J(R)\qq N_\F^I (T) = N_{\F'}^J (R)
\eqno £3.8.1;$$
hence, in order to prove our claim, it suffices to prove that $T$ is fully $I\-$nor-malized in $\F\,.$

\medskip
£3.9. We actually follow the proof of [5, Lemma~2.17];   for any $\F\-$morphism $\psi\,\colon T\.N_P^I
(T)\to P\,,$ set $Q' = \psi (Q)\,,$ denote by $\psi^*\,\colon Q'\cong Q$  the inverse
of the isomorphism $Q\cong Q'$ determined  by $\psi\,,$ and consider the $\F\-$morphism 
$$\iota_Q^P\circ \psi^* : Q'\too P
\eqno £3.9.1\phantom{.}$$
where  $\iota_Q^P\,\colon Q\to P$ is the inclusion map; it follows from [5,~Proposition~2.7] that
there is $\xi\,\colon N_P(Q')\to P$ such that $Q'' = \xi (Q')$ is both fully centralized and fully normalized in~$\F\,,$ and therefore,  since $Q$ is both fully centralized and fully normalized in $\F\,,$ it follows
from our argument above applied to $Q''$ and to $Q$ that there is an $\F\-$morphism 
$$\zeta : N_P (Q') \too P
\eqno £3.9.2\phantom{.}$$
mapping $Q'$ onto $Q\,.$ In particular, we have
$\zeta\big(N_P (Q')\big)\subset P'$ and, since $\psi \big(T\.N_P^I (T)\big)$ normalizes $Q'\,,$
 the homomorphism 
$$\eta : T\.N_P^I (T) = Q\.\big(R\.N_{P'}^J (R)\big)\too P'
\eqno £3.9.3\phantom{.}$$
mapping $w\in T\.N_P^I (T)$ on $\zeta \big(\psi (w)\big)$ belongs to $\F\big(P',T\.N_P^I (T)\big)\,;$
moreover, since $R\.N_{P'}^J (R)\subset  T\.N_P^I (T)$ and $\zeta \big(\psi (Q)\big) =Q\,,$ $\eta$ determines an
$\F'\-$morphism from $N_{P'}^J (R)$ to~$P'$~(cf.~£1.12.1); hence, since $R$ is fully $J\-$normalized in $\F'\,,$ 
we get (cf.~£3.8.1)
$$\zeta\Big(\psi \big(N_P^I (T)\big)\Big) =\eta \big(N_{P'}^J (R)\big) = 
N_{P'}^{{{}^{\eta}\!}J} \big(\eta (R)\big)\j \zeta
\Big(N_P^{{{}^{\psi\!}}I} \big(\psi (T)\big)\Big)
\eqno £3.9.4\phantom{.}$$
which forces $\psi \big(N_P^I (T)\big) = N_P^{{{}^{\psi}\!}I} \big(\psi
(T)\big)\,,$ proving the claim.

\medskip
£3.10. Consequently, if $\F'\not= \F$ then it follows from the induction hypothesis that $\F'$ is a Frobenius
$P'\-$category and, in particular, it fulfllls the corresponding condition~£2.3.2; thus, since $Q$ is still fully 
normalized and fully centralized in $\F'\,,$ for any subgroup $R$ of $P' = N_P(Q)$ containing $Q\.C_{P'}(Q)\,,$
the restriction induces a surjective group homomorphism
$$\F (R)_Q = \F' (R)_Q\too N_{\F' (Q)}\big(\F_R (Q)\big) = N_{\F (Q)}\big(\F_R (Q)\big)
\eqno £3.10.1,$$
so that, in this case, $\F$ also fulfills condition~£2.3.2.

\medskip
£3.11. Finally, assume that $P' = P$ and $\F' =Ê\F\,;$ we claim  that any {\it $\F\-$intersected\/} subgroup $R$ 
of $P$ (cf.~£1.11) contains $Q\,;$ indeed, since we have 
$$\F (P, R) = \big(N_\F (Q)\big) (P,R)
\eqno £3.11.1,$$
 any $\psi\in \F(P,R)$ can be extended to some $\hat\psi\in \F (P,Q\.R)$ and therefore we have  
$$\hat\psi\big(N_Q(R)\big)\i N_P \big(\psi (R)\big)
\eqno £3.11.1,$$
so that we get $\F_Q(R)\i {}^{\psi^*} \!\F_P\big(\psi (R)\big)\,;$ thus, according to equality~£1.11.2, 
we still have $\F_Q (R)\i \F_R (R)$ and therefore $N_Q (R)\i R\,,$ which forces $Q\i R\,.$

\medskip
£3.12. Firstly assume that $Q$ is {\it not  $\F\-$intersected\/}; then, we claim that $\F$ fulfills the two conditions
 in~Theorem~£2.7 above, so that $\F$ is a Frobenius $P\-$category  and, in particular, it fulfills condition~2.3.2.
According to~£3.11 and to  our choice of $Q\,,$ any $\F\-$intersected subgroup of $P$ belongs to $\frak X$ and therefore, since we are assuming that $\F^{^\frak X}$
is a partial Frobenius $P\-$category,  condition~£2.7.1 holds.
\eject

\medskip
£3.13.  Moreover, since any $\F\-$essential subgroup $U$ of $P$ is $\F\-$intersected (cf.~£1.11), $U$ belongs to $\frak X$ and we claim that any divisible $P\-$category $\hat\F$ fulfilling $\hat\F (P,U)\j\F (P,U)$
for every $\F\-$essential subgroup $U$ of $P$ contains~$\F\,;$ indeed, let $R$ be a subgroup of $P$ and
$\psi\,\colon R\to P$ an $\F\-$morphism; we may assume that $R$ is not $\F\-$essential
and then, as in~£3.4 above, it follows from Lemma~£1.5, Proposition~£1.9 and the Alperin condition that we have
$$\eqalign{\iota_R^P &= \iota_{U_0}^P \circ \nu_0\cr
\iota_{U_{i-1}}^P\circ \sigma_{i-1} \circ \nu_{i-1} 
&= \iota_{U_{i}}^P\circ \nu_{i}  \;\; \hbox{for any $1\le i\le \ell$}\cr 
\iota_{U_{\ell}}^P\circ \sigma_{\ell} \circ \nu_{\ell}   &= \psi\,.\cr}
\eqno £3.13.1\phantom{.}$$ 
for some $\ell$ and, for any $i\in \Delta_\ell\,,$ a suitable $\F\-$essential subgroup $U_i$ of $P$ and some 
elements $\sigma_i\in \F (U_i)$ and $\nu_i\in \F (U_i,R)\,;$ then, since $\hat\F$ is divisible, we have $\nu_0 = \iota_R^{U_0}$ and, in particular, it belongs to $\hat\F(U_0,R)\,;$ arguing by induction on $\ell\,,$ 
we may assume that $\nu_{\ell-1}\in  \hat\F (U_{\ell-1}, R)$ and, since  $\sigma_{\ell-1}$ belongs to~$\F(U_{\ell-1})\i \hat\F (U_{\ell-1})\,,$ we get $\nu_\ell
\in \hat\F (U_{\ell}, R)\,,$ 
so that $\psi$ belongs to~$\hat\F (P,R)$ since $\sigma_{\ell}\in \F(U_{\ell})\i \hat\F (U_{\ell})\,.$

\medskip
£3.14. Secondly,  assume that $Q$ is {\it $\F\-$intersected\/}; since we are assuming that $N_\F (Q) = \F\,,$
it is easily checked that, in this case, equality~£1.11.2 forces ${\bf O}_p \big(\F (Q)\big) = \F_Q (Q)\,;$
moreover, since $Q$ is $\F\-$selfcentralizing, in order to prove that condition~2.3.2 holds it suffices to consider 
a subgroup $R$ of~$P$ strictly containing $Q$ and then we have $\F_R (Q)\not= \F_Q (Q)\,,$ so that the
normalizer $K = N_{\F (Q)}\big(\F_R (Q)\big)$ is a proper subgroup of $\F(Q)\,.$

\medskip
£3.15. At present, set $P'' = N_P^K (Q)$ and $\F'' = N_\F^K (Q)\,;$ note that $P''$ contains $R\,,$
that $\F''$ fulfills the Sylow and the Alperin conditions (cf.~£3.2), and that we have $\F'' (Q) = K\,;$ since we also have
$$\F_Q (Q)\not= \F_R (Q) \triangleleft \F'' (Q) = K
\eqno £3.15.1,$$
$Q$ is {\it not\/} $\F''\-$essential; thus, any nonidentity element $\varphi\in \F'' (Q)$ defines a 
{\it $\F''\-$reducible $\F''\-$dimorphism\/} $\iota_Q^{P''}\circ( \varphi - {\rm id}_Q)$ and therefore, as in~£3.4 above, it follows from Lemma~£1.5, Proposition~£1.9 and the Alperin condition that we have
$$\eqalign{\iota_Q^{P''} &= \iota_{U_0}^{P''} \circ \nu_0\cr
\iota_{U_{i-1}}^{P''}\circ \sigma_{i-1} \circ \nu_{i-1} 
&= \iota_{U_{i}}^{P''}\circ \nu_{i}  \;\; \hbox{for any $1\le i\le \ell$}\cr 
\iota_{U_{\ell}}^{P''}\circ \sigma_{\ell} \circ \nu_{\ell}   &= \iota_Q^{P''}\circ \varphi\,.\cr}
\eqno £3.15.2\phantom{.}$$ 
for some $\ell$ and, for any $i\in \Delta_\ell\,,$ for a suitable $\F''\-$essential subgroup $U_i$ of $P''$ and
some elements $\sigma_i\in \F'' (U_i)$ and $\nu_i\in \F'' (U_i,Q)\,;$ note that, since we have $U_i\i P''\,,$ the image of $U_i$ in $\F (Q)$
normalizes $\F_{\!R} (Q)$ and therefore, since $Q$ is $\F\-$selfcentralizing, $U_i$~normalizes $R\,.$

\medskip
£3.16.  Then, for any $i\in\Delta_\ell\,,$ we claim that we can apply Lemma~£2.5
to $\F^{^\frak X}\!$ and to the triple $(Q,\F_R (Q),\varphi_i)$ where $\varphi_i\,\colon Q\to P$ is the $\F\-$morphism
defined by the restriction of~$\sigma_i\,;$ indeed, $Q$ is a normal proper subgroup of~$U_i\,,$ $U_i$ stabilizes
$\F_R (Q)$ and the $\F\-$morphism $\iota_{U_i}^P\circ \sigma_i\,\colon U_i\to P$ extends $\varphi_i\,.$
Consequently, it follows from this lemma that this triple is {\it extensile\/} and therefore,
 since $N_P^{\F_R(Q)}(Q) = R\,,$ there exists an $\F\-$morphism $\psi_i\,\colon R\to P$ extending~$\varphi_i\,;$
 moreover, since $\varphi_i$ is the restricion of $\sigma_i\in \F'' (U_i)\,,$ $\varphi_i$ normalizes $\F_R(Q)$ and therefore, since $Q$ is $\F\-$selfcentralizing,  we get $\psi_i (R) =R\,.$ Finally, the composition of the family of $\F\-$automorphisms of $R$ determined by 
 $\{\psi_i\}_{i \in \Delta_\ell}$ coincides with $\varphi\,;$ that is to say, the group homomorphism
 $$\F (R)\too N_{\F (Q)} \big(\F_R (Q)\big)
 \eqno £3.16.1\phantom{.}$$
 induced by the restriction is surjective, proving condition~£2.3.2. We are done.

\bigskip
\noindent
{\bf References}
\bigskip
\noindent
\smallskip\noindent
[1]\phantom{.} Jon Alperin,  {\it Sylow intersections and fusion\/}, Journal of
Algebra, 1(1964), 110-113.
\smallskip\noindent
[2]\phantom{.} James Green, {\it Functors on categories of finite group
representations\/}, Journal of Pure and Applied Algebra, 37(1985), 265-298.
\smallskip\noindent
[3]\phantom{.} Llu\'\i s Puig, {\it Structure locale dans les groupes
finis\/}, Bull.~Soc.~Math.~France, M\'emoire N$^o\,$47(1976)
\smallskip\noindent
[4]\phantom{.} Llu\'\i s Puig, {\it Source algebras of $p\-$central group 
extensions\/}, Journal of Algebra, 235(2001), 359-398.
\smallskip\noindent
[5]\phantom{.} Llu\'\i s Puig, {\it ``Frobenius categories versus Brauer blocks''\/}, Progress in Math.
274(2009), Birkh\"auser, Basel.

\end